\documentclass[a4paper,11pt]{article}
\usepackage[english]{babel}
\usepackage{epsf}
\usepackage[dvips]{epsfig}
 \usepackage{graphicx}
\usepackage{amssymb}
\usepackage{amsmath}
\usepackage[nice]{nicefrac}
\usepackage{amsfonts}
\usepackage{dsfont}
\usepackage[T1]{fontenc}
\usepackage{wasysym}
\usepackage{amssymb}
\usepackage{stmaryrd}
\usepackage{aeguill}

\newtheorem{theoreme-anglais}{Theorem}[section]

\newtheorem{cor-anglais}[theoreme]{Corollary}

\newtheorem{lemme-anglais}[theoreme]{Lemma}

\renewcommand{\epsilon}{\varepsilon}
\def\btab{\begin{eqnarray*}}
\def\etab{\end{eqnarray*}}
\def\beq{\begin{equation}}
\def\eeq{\end{equation}}

 \newcommand{ \un }{\mathds{1}}
 \newcommand{ \p }{\mathbb{P} }
 \newcommand{ \tpa }{\p^{\mathcal{E}} }

 \newcommand{ \E }{\mathbb{E}}

 \newcommand{ \R }{ \mathbb{R} }
 
 \newcommand{\N}{ \mathbb{N} }

\newcommand{ \px }{\overset{\leftarrow }{x}}
\newcommand{ \pz }{\overset{\leftarrow }{z}}
\newcommand{ \qz }{\overset{\rightarrow }{z}}
\newcommand{ \Ue }{U_{\epsilon}}

 \newcommand{ \lo }{ \mathcal{L} }

\newcommand{ \pn }{ \p\ a.s. - \mathcal{N}}
\newcommand{ \pne }{ P\ a.s. - \mathcal{N}}

\newtheorem{The}{{\bf Theorem}}[section]
 \newtheorem{Lem}[The]{Lemma}
 
 \newtheorem{Pro}[The]{\bf Proposition}

 \newenvironment{Pre}{\noindent \textbf{Proof.} \\ }{$\
 \blacksquare$}

\makeatletter

\@addtoreset{equation}{section} \makeatother

\title{The number of generations entirely visited for recurrent random walks on random environment
\author{P. Andreoletti, P. Debs \footnote{Laboratoire MAPMO - C.N.R.S. UMR 7349 - F\'ed\'eration Denis-Poisson, Universit\'e d'Orl\'eans
(France). \newline \vspace{0.1cm}  $\quad$  MSC 2000  60J55 ;  60J80 ; 60G50 ; 60K37. \newline \vspace{0.5cm} \textit{Key words :  random walks, random environment, trees} }
}}

\begin{document}
\maketitle
\bibliographystyle{plain}

\begin{abstract}
 In this paper we are interested in a random walk in a random environment on a super-critical Galton-Watson tree. We focus on the recurrent cases already studied by Y. Hu and Z. Shi \cite{HuShi10}, \cite{HuShi10a},  G. Faraud, Y. Hu and Z. Shi \cite{HuShi10b}, and G. Faraud \cite{Faraud}. We prove that the largest generation entirely visited by these walks behaves like $\log n$ and that the constant of normalization which differs from a case to another is function of the inverse of the constant of Biggins'  law of large number for branching random walks \cite{Biggins}.
 \end{abstract}

\section{Introduction and results}

First, let us  define the process:  \\
\textit{The environment $\textbf{E}$:} Let $\mathbb{T}_0$ a $N_0$-ary regular tree rooted at $\phi$. For all vertices $x  \in \mathbb{T}_0$ we associate a random vector $(A(x^{1}),A(x^{2}), \cdots, A(x^{N_x}),N_x)$ where $N_x$ is a non-negative integer bounded by $N_0$. We assume that the sequence $(A(x^{1}),A(x^{2}), \cdots, A(x^{N_x}),N_x), x \in \mathbb{T}_0)$ is i.i.d. and   that each vector has the same law as $(A_1,A_2,\cdots,A_N,N)$, we also assume that all $A_i$'s are independent of $N$. The sub-tree ${\mathbb T}=\{x\in \mathbb{T}_0, N(x)\neq 0 \}$ is a Galton-Watson tree (GW), so $(x^{1},x^{2}, \cdots, x^{N_x})$, are the $N_x$ children of $x$, and we denote $|x|$ the generation of $x$. For all  vertex $x$ in  $\mathbb {T}$, we denote $\px$ the parent of $x$, we also assume that $\phi$ has a unique ancestor denoted $\overset{\leftarrow }{\phi}$.   The set of environments denoted $\textbf{E}$ is the set of all sequences  $(A(x^{1}),A(x^{2}), \cdots, A(x^{N_x}),N_x), x \in \mathbb{T}_0)$,  we denote by $P$ the associated probability measure, and by $E$ the expectation. \\\textit{A random walk on $\mathcal E \in \textbf{E}$:}  we define a nearest neighbors random walk $(X_n,n\in \N,X_0=\phi )$ by its transition probabilities, 
\begin{align*}
& p(x,x^{i})= A(x^{i})/\left (\sum_{j=1}^{N_x}A(x^{j})+1 \right),\ p(x,\px)=1-\sum_{i=1}^{N_x}p(x,x^{i}),  \\
& p(\overset{\leftarrow }{\phi}, \phi)=1,
\end{align*}
note also that if $N_x=0$, then $p(x,\px)=1$.
We denote by $\p^{\mathcal{E}}$ the probability measure associated to  this walk, the whole system is described under the probability $\p$ which is the semi-direct product of measures $P$ and $\p^{\mathcal{E}}$. \\
\textit{General properties for the environment:} Note that by construction the GW is locally bounded, and we also add an ellipticity condition on the $A_i$'s,  
\begin{align}
& P-a.s\  \exists \ 0<\epsilon_0<1, \ \forall i, \epsilon_0 \leq A_i \leq 1/ \epsilon_0, \label{hyp1}
\end{align}
so the  moment-generating function  $\psi$ we define now, which contains the characteristics of the environment, is defined for all $t$:
$$\psi(t)=\log E\left(\sum_{i=1}^N A_i^t\right).$$
These assumptions (for the $A_i$'s, $1/A_i$'s and $N$), may be weaken by assuming exponential moments for all of them instead of ellipticity, but we do not think that we could reach easily the even weaker assumptions like in \cite{HuShi10b} for example. Nevertheless, we keep more generalist proofs as often as possible.
As mentioned in the abstract we assume that $\psi(0)>0$ so our Galton-Watson is super-critical, also that the random environment is non-degenerate.

\begin{figure}[ht]
\begin{center}
\input{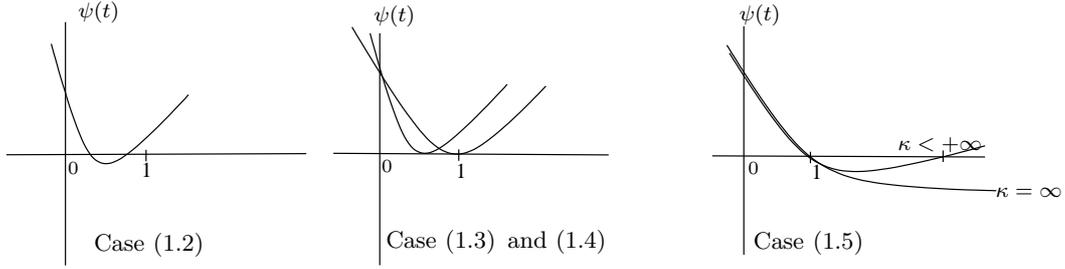} 
\caption{ The recurrent cases } \label{fig1}
\end{center}
\end{figure}

\noindent \textit{The recurrence criteria:} 
on a regular tree, they are first due to \cite{LyonPema}, in the present settings, we refer to (\cite{MenPet}) and the first part of \cite{Faraud}. Let 
$$\chi:=\inf_{t \in [0,1]}\psi(t),$$ 
 then the walk is transient if and only if  $\chi >0 $. The recurrent case can be specified as follows, if 
\begin{align}
\chi< 0 \label{hypP}
\end{align}
then the walk is positive recurrent, to determine the other case, we have to take into account the sign of  
$$ \psi'(1)=e^{-\psi(1)} E \left[\sum_{i=1}^NA_i \log A_i \right].$$
If
\begin{align}
\chi &= 0 \textrm{ and }   \psi'(1)>0, \label{hypPP}
\end{align}
the walk is positive recurrent, whereas if
\begin{align}
\chi &= 0 \textrm{ and }   \psi'(1)=0, \label{hyp0} \textrm{or} \\
\chi &= 0 \textrm{ and }   \psi'(1)<0,  \label{hyp00}
\end{align}
the walk is null recurrent.  In figure (\ref{fig1}) we present the shape of  $\psi$ for each case, for the last one (\ref{hyp00})  a constant appears naturally: 
 \begin{align*}
\kappa:=\inf\{t>1, \psi(t)=0 \} \in (1,+ \infty].
\end{align*} 
\textit{Asymptotics for the largest visited generation $X_n^*$:} The asymptotic behavior of $X_n^*:= \max_{1 \leq k \leq n} \left| X_k \right|$ is well known thanks to the works of Y. Hu and Z. Shi \cite{HuShi10a}, \cite{HuShi10b} and G. Faraud, Y. Hu and Z. Shi \cite{HuShi10}. They prove that there is three main different behaviors, the first one   (\cite{HuShi10a}) says that the walk is very slow and will never reach a generation larger than $\log n$ for an amount of time $n$, more pricisely
\begin{align*} 
& \textrm{ if } (\ref{hypP}) \textrm{ is realized then  } \pn,  \lim_{n \rightarrow +Ê\infty}\max_{ 0 \leq i \leq n}\frac{|X_i|}{\log n}=C_1, 
\end{align*}
where $\pn$ means $\p$ almost surely on the set of non-extinction of the Galton Watson tree.  Note that in \cite{HuShi10a} a regular tree is considered but the result remains true with our hypothesis. In \cite{HuShi10b} and \cite{HuShi10a}, it is proven that
\begin{align*} 
& \textrm{ if } (\ref{hypPP})  \textrm{ is realized then } \pn, \lim_{n \rightarrow +Ê\infty}\max_{ 0 \leq i \leq n}\frac{|X_i|}{(\log n)^3}=C_2, \\
& \textrm{ if } (\ref{hyp0})  \textrm{ is realized then } \pn, \lim_{n \rightarrow +Ê\infty} \max_{ 0 \leq i \leq n}\frac{|X_i|}{(\log n)^3}=C_3,
\end{align*}
in this delicate case, there is still a slow movement, but they prove that the environment allows enough regularity to let the walk escape until generation $(\log n)^3$. Note that in \cite{HuShi10} they work with a more general setting, a GW tree, weaker hypothesis of regularity than ours,  and succeed to determine $C_2$ and $C_3$.  Finally, there is also a sub-diffusive case  also obtained in (\cite{HuShi10a}) : 
\begin{align*} 
& \textrm{ if } (\ref{hyp00})  \textrm{ is realized then } \pn, \lim_{n \rightarrow +Ê\infty}  \frac{  \log \max_{ 0 \leq i \leq n} |X_i| } {\log n}= 1- \frac{1}{\min{(\kappa,2)}}. 
\end{align*}
Note also that for large $\kappa$ \cite{Faraud} shows the existence of a central limit theorem for this last case. \\

\noindent In this paper we are interested in the \textit{largest generation entirely visited by the walk}, more precisely we get the asymptotic behavior of
\begin{align*}
& R_n:=\sup\lbrace k\geq1,\forall  \vert z\vert= k,\lo(z,n)\geq1 \rbrace, 
\end{align*}
with $\lo$ the local time of $X$ defined by $\lo(z,n):= \sum_{k=1}^n \un_{X_k=z}$.

We also need the following constant of law of large number for branching random walks :
\begin{align}
\tilde{J}(a)&:=\inf_{t\geq 0}\{\psi(-t)-at\}, \label{J} \  
\tilde{\gamma}:=\sup\{a \in \R,\ \tilde J(a)>0\}, 
\end{align}
note that as $\chi \leq 0$, $\tilde{\gamma}>0$. Our main result shows that, contrary to $X_n^*$, there is essentially two cases:
\begin{The}  Assume (\ref{hyp1}), then if (\ref{hypP}) or (\ref{hypPP}) or (\ref{hyp0}) are realized, $\pn$
$$\lim_{n\rightarrow+\infty}\frac{ R_n}{\log n}=\frac{1}{\tilde \gamma},$$
otherwise  $\pn$
$$\lim_{n\rightarrow+\infty}\frac{ R_n}{\log n}=\frac{1}{\tilde \gamma\min(\kappa,2) }.$$
\end{The}
So the largest generation entirely visited is far smaller than the largest generation visited by these walks except for the slowest case (\ref{hypP}). In fact there is no difference between the  first three  cases (which are the slowest ones) and we see appear the characteristic constant $\kappa$ for the fourth one. In fact, if instead of stopping the walk at a deterministic time $n$ we stop it at $n$ return time to the root, we have no longer any difference. More precisely, for all $i\geq 1$ let $T_{\phi}^i:=\inf\{k>T_{\phi}^{i-1}, X_k=\phi\}$ the $i^{th}$ return time to $\phi$, with $T_{\phi}^0=0$ and denote $\tilde R_n:=R_{T_\phi^n}$ then
\begin{Pro} \label{Pro1} Assume (\ref{hyp1}), then $\pn$
$$\lim_{n\rightarrow+\infty}\frac{\tilde R_n}{\log n}=\frac{1}{\tilde \gamma}.$$
\end{Pro}
This last fact shows that  the difference for all the cases appears only in the behavior of the local time at the root $\phi$. In fact we only need the logarithm behavior of $\lo$ at $\phi$ , it is given by
\begin{Pro}\label{profin}
Assume (\ref{hyp1}), then if (\ref{hypP}) or (\ref{hypPP}) or (\ref{hyp0}) are realized, $\pn$
$$\lim_{n\rightarrow+\infty}\frac{ \log \lo(\phi,n)}{\log n}=1,$$
otherwise  $\pn$
$$\lim_{n\rightarrow+\infty}\frac{ \log \lo(\phi,n)}{\log n}=\frac{1}{\min(\kappa,2) }.$$
\end{Pro}

\noindent (\ref{hypP}) and (\ref{hypPP}) are obvious given recurrence positivity. 

The rest of the paper is organized as follows, in Section 2, we prove the result for $\tilde R_n$, it is the upper bound that needs more attention. In Section 3 we move from $\tilde R_n$ to $R_n$, also for the sake of completeness we add classical results in an appendix.

To study asymptotical behaviours associated to $(X_n)_{n\in\mathbb N}$, a quantity appears naturally: the potential process $V$ associated to the environment which is actually a branching random walk. It is defined by $V(\phi):=0$ and
$$V(x):=-\sum_{z\in\rrbracket \phi, x  \rrbracket } \log A(z),x\in\mathbb T\backslash \lbrace \phi\rbrace, $$
where $\llbracket \phi, x  \rrbracket  $ is the set of vertices on the shortest path connecting $\phi$ to $x$ and $\rrbracket \phi, x  \rrbracket=\llbracket \phi, x  \rrbracket\backslash \lbrace \phi\rbrace$.  

\section{Proof of Proposition \ref{Pro1}}

\subsection{Lower bound}

In this first section, we prove that $\pn$ for $n$ large enough
\begin{equation}\label{liminf}
\frac{\tilde R_n}{\log n}\geq \frac{1-\varepsilon}{\tilde \gamma}=:c_1.
\end{equation} 
For this purpose, note that:
\begin{eqnarray*}
\tpa_{\phi}(\tilde R_n< c_1\log n)=\tpa_{\phi}\left(\bigcup_{\vert z\vert =c_1\log n}\lbrace \lo(z,T_\phi^n)=0 \rbrace\right)=\tpa_{\phi}\left(A_n \right)
\end{eqnarray*}
where $A_n:=\bigcup_{\vert z\vert =c_1\log n}\lbrace T_z>T_\phi^n\rbrace$. Note that for typographical simplicity, we do not make any difference between a real number and its integer part. Thus, according to strong Markov property:
\begin{eqnarray*}
\tpa_{\phi}(A_n) \leq \sum_{\vert z\vert =c_1\log n}\tpa_{\phi} (T_z>T_\phi)^n \leq Z_{c_1\log n}\max_{\vert z\vert =c_1\log n}e^{n\log \tpa_{\phi}(T_z>T_\phi)}
\end{eqnarray*}
 with $Z_n:=\mbox{Card}\lbrace\vert z\vert= n  \rbrace$, the number of vertices in the $n$-th generation. With $\E[Z_1]=\E[N]=e^{\psi(0)}$, the expected number of offspring at the first generation, it is a classical result that $W_n:=\frac{Z_n}{e^{n\psi(0)}}$ is a positive martingale an consequently $(W_n)_{n\geq 0}$ admits a.s. a limit  when $n$ goes to infinity. So, there exists $C(\omega)$ and $n_0(\omega)$ such that: $\forall n\geq n_0(\omega),\, \frac{Z_n}{e^{n\psi(0)}}\leq C(\omega).$ 
Consequently  $\forall n\geq n_0(\omega)$, noting that ${e^{\psi(0)c_1\log n}=n^{c_1\psi(0)} }$:
\begin{equation}\label{rmaj}
\tpa_{\phi}(A_n)\leq  C(\omega) n^{c_1\psi(0)} \max_{\vert z\vert =c_1\log n}e^{n\log(1- \tpa_\phi(T_z<T_\phi))}.
\end{equation}
As $X$ is recurrent,  $\tpa_\phi(T_z<T_\phi)$ tends to 0 when $n$ goes to infinity 
and we have to study the asymptotical behaviour of: 
$$\aleph_n:=\max_{\vert z\vert =c_1\log n}e^{-n \tpa_\phi(T_z<T_\phi)}=\max_{\vert z\vert =c_1\log n}e^{-np(\phi,\phi_z) \tpa_{\phi_z}(T_z<T_\phi)},$$
where $\phi_z$ is the child of $\phi$ in $\rrbracket \phi,z\rrbracket $. \\
Recall that, thanks to the ellipticity conditions, $ \forall u\in \mathbb T,\, e^{-V(u)}=A(u)>\varepsilon_0$, formulas (\ref{eq1}) yields:
\begin{equation}\label{f1}
\tpa_{\phi_z}(T_z<T_\phi)=\frac{e^{V(\phi_z)}}{\sum_{u\in\rrbracket \phi,z \rrbracket}e^{V(u)}}\geq \varepsilon_0\frac{e^{-\overline{V}(z)}}{\vert z\vert }=\varepsilon_0\frac{e^{-\overline{V}(z)}}{c_1\log n},
\end{equation}
where $\overline{V}(z)=\max_{x\in\rrbracket \phi,z \rrbracket}V(x)$.  The ellipticity conditions ensure that there is a constant $K>0,$ such that   $\forall z\in \mathbb T,\,K<\nicefrac{\varepsilon_0 p(\phi,\phi_z)}{c_1}$, then using \ref{f1}:
\begin{equation}\label{finprop1}
\aleph_n\leq  \max_{\vert z\vert =c_1\log n}e^{-\frac{n\varepsilon_0p(\phi,\phi_z)}{c_1\log n}e^{-\overline{V}(z)}}\leq  e^{-\frac{Kn}{\log n}e^{-\max_{\vert z\vert =c_1\log n}\overline{V}(z)}} .
\end{equation} 

At this level, it remains to study $\overline V$ and we need the following:
\begin{Lem} \label{lemmax}Assume $\chi \leq 0$, there exists a constant $a>0$ such that $\pne$ for $\ell$ large enough : 
$$\max_{|z|=\ell} \overline{V}(z) \leq \tilde{\gamma} \ell \left(1+ a \frac{\log \ell}{\ell} \right) $$
\end{Lem}
Let us postpone the proof of this lemma and finish the proof of (\ref{liminf}): for $n$ large enough, the previous lemma implies: 
$$\max_{|z|=c_1\log n} \overline{V}(z) \leq \tilde{\gamma}  c_1\log n \left(1+ \epsilon/2 \right),$$
and one can write $\pne$ for $n$ large enough:
\begin{equation}\label{finprop2}
\aleph_n\leq e^{-\frac{Kn^{1-c_1\tilde\gamma(1+\frac{\varepsilon}{2})}}{\log n}}\leq  e^{-\frac{Kn^\frac{\varepsilon}{2}}{\log n}}. 
\end{equation}
Finally formulas (\ref{rmaj}) and (\ref{finprop2}) give that $P$ almost surely on the set of non-extinction $\sum\tpa_{\phi}(A_n)<\infty$, thus (\ref{liminf}) is established using Borel-Cantelli Lemma.  
\noindent \\

\noindent \textbf{Proof of lemma \ref{lemmax}:}\\
This result is classical and for the sake of completeness, we give some details below. Let $\epsilon_\ell:=a{\log \ell}/{\ell}$, 
using the Biggins identity (\ref{B1}), we easily obtain:
\begin{align*}
P\left(\max_{|z|=\ell} \overline{V}(z) \geq \tilde\gamma \ell (1+ \epsilon_\ell)\right) &= P\left(\cup_{j=1}^\ell\cup_{|z|=j}{\left\{V(z) \geq \tilde\gamma \ell(1+Ê\epsilon_\ell)\right\}}\right) \\
\leq  \sum_{j=1}^\ell E\left(\sum_{|z|=j}\mathds{1}_{\left\{V(z) \geq \tilde\gamma \ell(1+Ê\epsilon_\ell)\right\}}\right) 
&=  \sum_{j=1}^\ell e^{j \psi(1)}E\left(e^{S_j}\un_{\left\{S_j \geq \tilde\gamma \ell(1+Ê\epsilon_\ell)\right\}}\right).
\end{align*}
For any  $b>0$,  a simple partition of the event $\left\{S_j \geq \tilde\gamma \ell(1+Ê\epsilon_\ell)\right\}$ gives: 
\begin{align*}
 E\left[e^{S_j}\un_{\left\{S_j \geq \gamma \ell(1+Ê\epsilon_\ell)\right\}}\right] &= \sum_{r=0}^{+ \infty}E\left[e^{S_j}\un_{\left\{S_j \in [\tilde\gamma \ell(1+Ê\epsilon_\ell)+br,\tilde\gamma \ell(1+Ê\epsilon_\ell)+b(r+1)[\right\}}\right] \\
  & \leq \sum_{r=0}^{+ \infty}e^{\tilde\gamma \ell(1+Ê\epsilon_\ell)+b(r+1)} \p\left(S_j \geq \tilde\gamma \ell(1+Ê\epsilon_\ell)+br \right).
\end{align*} 
The ellipticity condition gives $e^{\psi(-\delta)}=E[\sum_{\vert x\vert=1}e^{ \delta V(x)}]\leq \left(\frac{1}{\epsilon_0}\right)^\delta E[N]<\infty$ for all $\delta \in \R$, so according Biggins identity (\ref{B2}), $\E[e^{(1+\delta)S_1}]<+ \infty$. Thus, using Markov inequality and the fact that $(S_{i}-S_{i-1},i\geq1)$ are i.i.d. random variables, 
$\forall c>0,\,\p(S_j\geq c)\leq  \frac{\E[e^{(1+ \delta) S_1}]^j}{e^{(1+ \delta) c}} .$
Collecting the previous inequalities, and taking $c= \tilde \gamma \ell(1+Ê\epsilon_\ell)+br$:
\begin{eqnarray}
\p\left(\max_{|z|=\ell} \overline{V}(z) \geq \tilde \gamma \ell (1+ \epsilon_\ell)\right)&  \leq & e^{b-\delta \tilde \gamma \ell(1+\varepsilon_\ell)}\sum_{ r \geq0}e^{- \delta r b}\sum_{j=1}^\ell\E[e^{(1+\delta)S_1}]^je^{\psi(1) j} \nonumber \\
&=& \frac{e^b }{1-e^{-\delta b}}  e^{-\delta \tilde \gamma \ell(1+\varepsilon_\ell)}  \sum_{j=1}^\ell e^{\psi(-\delta) j } \nonumber \\
&=& \frac{e^b}{1-e^{-\delta b}} \frac{ e^{\psi(-\delta)}}{e^{\psi(-\delta)}-1} e^{-\delta \tilde \gamma \ell(1+\varepsilon_\ell)}{\left(e^{\psi(-\delta) \ell }-1\right)} \nonumber \\
&\leq & M  e^{-\delta \tilde \gamma \ell \varepsilon_\ell }{e^{\ell(\psi(-\delta)-\delta \tilde \gamma)  }} =: M\Delta_\ell(\delta),   \label{eq000} 
\end{eqnarray}
for the first equality we use Biggins identity, for the second one the fact that for all $\delta>0$, $e^{\psi(-\delta)}> e^{\psi(0)}>1 $ and $M$ is a positive constant.\\
 Before going any further, according to the definitions of $\tilde{J}$  and $\tilde{\gamma}$ see (\ref{J}), note that $\tilde J(\tilde \gamma)=0$. Indeed $\psi$, as a function of $t$, is convex moreover by hypothesis $\psi(0)>0$ and $\inf_{t \in [0,1]} \psi(t) \leq 0$, so it reaches its minimum for some $t>0$, so $\tilde{J}(0)=\inf_{t \geq 0} \psi(-t)>0$. Moreover by hypothesis $\psi(-t)$ is finite for every $t>0$, and therefore for all $t$ we can find some $a$, large enough such that $-\infty<\tilde J(a)\leq\psi(-t)-ta<0$. Then the definition of $\tilde \gamma$ gives effectively that  $\tilde J(\tilde \gamma)=0$.
We can now come back to $\Delta_\ell$, we have two cases, either
\begin{itemize}
\item there exists $t_0>0$ such that $ \psi(-t_0)-t_0\tilde \gamma=0$. Then $\sum_{\ell\geq0}\Delta_\ell(t_0)=M \sum_{\ell\geq0}e^{-t_0 \tilde\gamma \ell\epsilon_\ell}<\infty$, and we conclude with the Borel-Cantelli Lemma, or
\item $ \psi(-t)\sim\tilde \gamma t$ when $t$ goes to infinity, note that by convexity of $\psi$, $\psi(t)-\tilde \gamma t \geq 0$ for all $t$. Then we can take $\delta=\delta_\ell=\frac{1}{\varepsilon_\ell}$, in this case $\Delta_\ell(\delta_\ell)\sim e^{-\ell\tilde\gamma}$ and we easily conclude with Borel-Cantelli Lemma.\hfill $\blacksquare$
\end{itemize}

\subsection{ Upper bound}
In this section we prove that, for all $\epsilon>0$, $\pn$ for all $n$ large enough
\begin{align}
\frac{\tilde{R}_{n}}{ \log n} \leq  c_2:=\frac{1+\varepsilon}{\tilde \gamma}.  \label{limsup}
\end{align}
The strategy is the following, we first make a first cut in the tree close to the root at a generation which depends on $\epsilon$. We denote $(z_i,i \leq \Ue)$ the vertices of this generation of the tree. We show that during the $n$ return time to the root the local time at each of these individuals is not much larger than $n$ (Lemma \ref{lem2.2}). Then we make a second cut in the tree at generation $(1+ \epsilon/2) \log n$. We select at this generation one descendant for each $z_i$ called $\underline{z}_i$ satisfying the property to have a large potential $V(\underline{z}_i)$ (see \ref{alakokan}). We prove that the local times on these vertices during the return time to $z_i$ do not exceed a power of $\log n$ almost surely (Lemma \ref{2.3}). We finally prove a last technical lemma (Lemma \ref{lem2.5}) which shows that there are very few back and forth movements between $z_i$ and its descendant $\underline{z}_i$.  
Finally,  using the three Lemmata we can extract some parts of the trajectory of the random walk (before the  $n$th visit to the root) which are independent up to a translation in time. Using this independence we finally prove that $ \tpa\left(\frac{\tilde{R}_{n}}{\log n} > c_2\right)$ is summable which leads to the result.
\\

Let $u_{\epsilon}$ a positive integer that will be precised later. Let $(z_i,i \leq \Ue=:|Z_{u_{\epsilon}}|)$, the individuals of generation $u_{\epsilon}$. We first prove that before the  $n$th visit to the root each point at generation $u_{\epsilon}$ can not be visited many more times than $n $. 

\begin{Lem} \label{lem2.2} Assuming (\ref{hypP}), for all positive and increasing sequence of integers  $(h_n,n \in \N)$ with $\lim_{n \rightarrow + \infty}h_n=+ \infty$, $\pne$ for $n$ large enough
\begin{align*}
\p_\phi^{\mathcal E}\left( \bigcup_{1 \leq j \leq \Ue} \left\{  \lo({z}_j,T_{\phi}^{n})\geq h_n  n \right\}  \right) \leq h_n2^{-n}.
\end{align*}
\end{Lem}

\begin{Pre}
Let us denote $ \bigcup_{1 \leq j \leq \Ue} \bar{\mathcal {A}}_j$ the event in the previous probability. Let $q_{z_j}>0$ and $r_{z_j}>0$ two sequences that we define later. Using successively Markov inequality and the strong Markov property: 
\begin{eqnarray}
 \p_\phi^{\mathcal E}(\lo(z_j,T_\phi^n)\geq r_{z_j}n) \leq e^{-q_{z_j} r_{z_j}  n}\E_\phi^{\mathcal E}\left[e^{q_{z_j}\lo(z_j,T_\phi^n)}\right]=e^{-q_{z_j} r_{z_j}  n} \left(\E_\phi^{\mathcal E}\left[e^{q_{z_j}\lo(z_j,T_\phi)}\right]\right)^n. \label{eqqp1}
\end{eqnarray}
Let us denote $w_{z_j}:=\p^{\mathcal E}_{z_j}(T_{z_j}>T_\phi)$, $v_{z_j}:=\p^{\mathcal E}_\phi(T_\phi>T_{z_j})$. Assuming that  for all $j,\ e^{q_{z_j}}(1-w_{z_j})<1 $:
\begin{eqnarray*}
\E_\phi^V\left[e^{q_{z_j}\lo(z,T_\phi)}\right]
=1-v_{z_j}+v_{z_j}e^{q_{z_j}}\frac{w_{z_j}}{1-(1-w_{z_j})e^{q_{z_j}}}=1+v_{z_j}\frac{e^{q_{z_j}}-1}{1-(1-w_{z_j})e^{q_{z_j}}}.
\end{eqnarray*}
As for all $j$, $1-w_{z_j}<1$ we can chose $q_{z_j}=\log(1+w_{z_j})$ which obviously satisfied  $e^{q_{z_j}}(1-w_{z_j})<1 $, we obtain:
\begin{eqnarray}
\E_\phi^{\mathcal E}\left[e^{q_{z_j}\lo(z,T_\phi)}\right]&=& 1+ \frac{v_{z_j}}{w_{z_j}}. \label{eqqp2} \end{eqnarray}
Replacing this expression in (\ref{eqqp1}), as $v_{z_j}\leq 1$: 
\begin{eqnarray}
 \p_\phi^{\mathcal E}(\lo(z_j,T_\phi^n)\geq r_{z_j}n) \leq  \left(\left(\frac{1}{1+w_{z_j}}\right)^{r_{z_j}}Ê\left(1+\frac{1}{w_{z_j}}\right)\right)^n,
\end{eqnarray}
finally taking $r_{z_j}=2 \log(1+1/w_{z_j})/\log(1+w_{z_j}) $, we get
\begin{align*}
\p_\phi^{\mathcal E}\left( \bigcup_{1 \leq j \leq \Ue} \left\{  \lo(z_j,T_\phi^n)\geq r_{z_j} n  \right\}  \right) & \leq \Ue \max_{j \leq \Ue} \p_\phi^{\mathcal E}(\lo(z_j,T_\phi^n)\geq r_{z_j}n) \leq   \Ue 2^{-n}.
\end{align*}
To finish, we have to estimate $r_{z_j}$ and so $w_{z_j}=p(z_j,\pz_j)\p^{\mathcal E}_{\pz_j}(T_{z_j} \geq T_\phi)$. By (\ref{eq2}) we note that $w_{z_j}$ can be small if the potential from the root to $z_j$ decreases, but thanks to the hypothesis of ellipticity, $P$ a.s.  $w_{z_j} \geq c' (\epsilon_0)^{\Ue}$, where $c'>0$ so  $P$ a.s.  $r_{z_j} \leq c'' (\epsilon_0)^{-2\Ue} $ with $c''>0$.  By the ellipticity condition for $N$, $P$ a.s. for $n$ large enough $r_{z_j} \leq h_n$, and $\Ue \leq h_n$, so $ \p_\phi^{\mathcal E}\left(\bigcup_{1 \leq j \leq \Ue} \bar{\mathcal {A}}_j \right)\leq \p_\phi^{\mathcal E}\left( \bigcup_{1 \leq j \leq \Ue} \left\{ \lo(z_j,T_\phi^n)\geq r_{z_j} n  \right\}  \right)\leq h_n2^{-n}$.
\end{Pre}
\\

In what follows, for simplicity, we denote $z>x$ if $\rrbracket x,z\rrbracket\neq \emptyset$, in other words $x$  is an ancestor of $z$.
\noindent \\ Let $(\underline{z}_i,i \leq \Ue)$ the individuals of generation $a_n:=(1+\epsilon/2)\log n/\tilde \gamma$  such that  $z_i<\underline{z}_i$ and satisfying that for all $1 \leq i \leq \Ue $:
\begin{align}
 V(\underline{z}_i)-V(\qz_i)\geq \tilde \gamma a_n \left(1-b\frac{\log a_n}{a_n}\right),\ \max_{u \in \rrbracket z_j,\underline z_j \rrbracket} {V}(u)-V({z}_j) \leq  \tilde \gamma c \log a_n,  \label{alakokan}
\end{align}
where ${\qz}_j$ the descendant of $z_j$ on $\rrbracket z_j,\underline{z}_j\rrbracket$. We prove in Lemma \ref{LemVL} below that such points exists almost surely. Define also 
$$ K_{n}={(\log n)^{3+c \tilde \gamma}}/{h_n n}. $$ 
We now prove that the probability for the local time, at each points $\underline{z}_i$ until  $T_{z_j}^{\lo({z}_j,T_{\phi}^{n})} $, to be larger than $n K_n$ is rather small.
 
\begin{Lem} \label{2.3} Assuming (\ref{hypP}), there exists  a constant $c_3>0$ such that $\pne$ for $n$ large enough
\begin{align*}
\p_\phi^{\mathcal E}\left( \bigcup_{1 \leq j \leq \Ue} \left\{  \lo\left(\underline{z}_j,T_{z_j}^{\lo({z}_j,T_{\phi}^{n})}\right)\geq K_n h_n n \right\}  \right) \leq  e^{- \frac{c_3}{4} (\log n)^2 }.
\end{align*}
\end{Lem}

\begin{Pre}
Let us denote $ \bigcup_{1 \leq j \leq \Ue} \bar{\mathcal B}_j$ the event in the previous probability, from  (\ref{eqqp1}) and (\ref{eqqp2}):
\begin{eqnarray*}
A:=\p_\phi^{\mathcal E}\left(\bigcup_{1 \leq j \leq \Ue} \left \{ \bar{\mathcal B}_j, \mathcal A_j \right \} \right)&\leq& \Ue \max_{ 1 \leq j \leq \Ue } \p_\phi^{\mathcal E}(\lo(\underline{z}_j,T_{z_j}^{h_n n})\geq K_{n} h_n n), \\
& \leq & \Ue \max_{ 1 \leq j \leq \Ue } \left(\left(\frac{1}{1+\tilde w_{z_j}}\right)^{K_{n}} \left(1+\frac{\tilde v_{z_j}}{\tilde w_{z_j}}\right)\right)^{nh_n},
\end{eqnarray*}
where $\tilde{w}_{z_j}:=\p^{\mathcal E}_{\underline z_j}(T_{z_j}< T_{\underline{z}_j})$ and $\tilde v_{z_j}:=\p^{\mathcal E}_{z_j}(T_{z_j}>T_{\underline{z}_j})$. Using Lemma \ref{lemfinal} and  the hypothesis of ellipticity, $P$ a.s.
\begin{align*}
 \tilde{v}_{z_j}  \leq e^{-(\max_{u \in \rrbracket z_j,\underline z_j \rrbracket} {V}(u)-V({\qz}_j))}, 
  \tilde{w}_{z_j}  \geq \frac{c_0'}{a_n}  e^{-(\max_{u \in \rrbracket z_j,\underline z_j \rrbracket} {V}(u)-V(\underline{z}_j))},
\end{align*}
with $c_0'>0$. \noindent Note that for all $0<c'<1$, and $x$ small enough $(1+x)^{-\alpha} \leq (1-c'\alpha x) $, taking $c'=1/2$, $x= \tilde w _j$, and $\alpha=K_n$, we get for all $n$ large enough: 
\begin{multline*}
A\leq \Ue  \max_{ 1 \leq j \leq \Ue }\left(\left(1-\frac{ \tilde w_{z_j}}{2} K_{n}\right) \left(1+\frac{ \tilde v_{z_j}}{\tilde w_{z_j}} \right)\right)^{h_n n} 
 \leq \Ue  \max_{ 1 \leq j \leq \Ue }\left(1-\frac{\tilde w_{z_j}}{2} K_{n}+\frac{\tilde v_{z_j}}{\tilde w_{z_j}} \right)^{h_n n} \\
 \leq  \Ue  \max_{ 1 \leq j \leq \Ue }\left(1- \frac{c_0^\prime K_{n}}{2a_n}{e^{-(\max_{u \in \rrbracket z_j,\underline z_j \rrbracket} {V}(u)-{V}(\underline{z}_j))}}+\frac{a_n}{c_0'} e^{-(V(\underline{z}_j)-V(\qz_j))} \right)^{h_n n}.
\end{multline*}
Now, assume for the moment that the sequence $(\underline{z}_j, j \leq \Ue)$ we have defined in (\ref{alakokan}) exists $\pne$, then $\pne$ for $n$ large enough 
\begin{eqnarray*}
A&\leq&\Ue \left(1-\frac{c_0'}{2}\frac{e^{-\tilde\gamma c \log a_n }}{\log n}K_n+\frac{a_n}{c_0'}  e^{-\tilde{\gamma}a_n+\tilde{\gamma}b \log a_n)}    \right)^{h_n n} \\
& \leq & \Ue \left (1-\frac{c_0'}{2}\frac{(\log n)^{2}}{ nh_n}+\frac{1}{c_0'}\frac{a_n^{1+\tilde \gamma b}}{n^{(1+\varepsilon/2)}}    \right)^{h_n n}  \leq \Ue e^{- \frac{c_0^{\prime\prime}}{2} (\log n)^2 }.
\end{eqnarray*}
To finish just notice that $\p_\phi^{\mathcal E}\left(\bigcup_{1 \leq j \leq \Ue}  \bar{\mathcal B}_j  \right) \leq \p_\phi^{\mathcal E}\left(\bigcup_{1 \leq j \leq \Ue} \left \{ \bar{\mathcal B}_j, \mathcal A_j \right \} \right)+ \p_\phi^{\mathcal E}\left(\bigcup_{1 \leq j \leq \Ue} \bar{\mathcal A}_j \right) $, use Lemma \ref{lem2.2} and the ellipticity condition for $N$.
\end{Pre} 

\noindent \\ We are left to prove the following  
\begin{Lem} \label{LemVL} Assume \ref{hyp00} then there exist two constants $b_0>0$ and $c_0>0$ such that $P$ almost surely on the set of non extinction, for all $l$ large enough 
\begin{eqnarray}
\exists z, |z|=\ell,\ V(z)\geq \tilde \gamma \ell\left(1-b_0\frac{\log \ell}{\ell}\right),\ \overline{V}(z)-V(z) \leq  \tilde \gamma c_0 \log \ell,  
\end{eqnarray}
For all integer $A>0$, let us denote $(z_i, 1 \leq i \leq U_{A})$ the individuals of the $A^{nth}$ generation, then there exist two constants $b>0$ and $c>0$ such that $\pne$ for all $l$ large enough  and all $i \leq U_{A}$ 
\begin{eqnarray}
\exists \underline{z}_i, |\underline{z}_i|=\ell,\ V(\underline{z}_i)-V(\qz_i)\geq \tilde \gamma \ell\left(1-b\frac{\log \ell}{\ell}\right),\ \max_{u \in \rrbracket z_j,\underline z_j \rrbracket} {V}(u)-V(\underline{z}_i) \leq  \tilde \gamma c \log \ell.  
\end{eqnarray}
\end{Lem}
\textbf{Proof}
First note that the second part of the Lemma is a simple consequence of the first part of the Lemma the ellipticity condition and the stationarity of the potential $V$.
So if we prove that there exists two constants $a>0$ and $b>0$ $P$ almost surely on the set of non extinction for $n$ large enough:
$$ \left\{\max_{|z|=\ell}\overline{V}(z) \leq  \tilde \gamma \ell(1+a \log \ell / \ell );\ \exists z,\ |z|=\ell,\ V(z) \geq \tilde \gamma \ell (1-b \log \ell /\ell) \right \}$$
then we get the first part of the Lemma. 
We have already proven, in Lemma \ref{lemmax}, that there exists a constant $a>0$ such that $P$-a.s. on the set of non extinction for $\ell$ large enough $\max_{|z|=\ell} \overline{V}(z) \leq \tilde \gamma \ell (1+a \log \ell / \ell )$. So we just need that  $P$-a.s. on the set of non extinction for $\ell$ large enough $ \exists z,\ |z|=\ell,\ V(z) \geq \tilde \gamma \ell (1-b \log \ell / \ell)$, for this we use the results of \cite{McDiarmid}, note that here we are interested in the maximum instead of the minimum so few changes occur. Let $\tilde F(t):=E\left[\sum_{|x|=1}\un_{V(x) \geq t}\right]$, by independence of $N$ and the increments $A_i$, we have $\tilde F(t)=\sum_{j=1}^{+ \infty}\sum_{i=1}^jP(N=j)P(-\log A_i \geq t)$ and by hypothesis (\ref{hyp1}), for all $t \geq - \log( \epsilon_0)$, $\tilde F(t)=0 $, therefore $\tilde \alpha:=\sup\{t, \ \tilde F(t)>0\} $ is finite. In \cite{McDiarmid} there is two theorems the first one and the remarks that follow concern the case with a finite $\tilde \alpha$ and $\tilde F( \tilde \alpha )\geq 1$ and the second one the case $\tilde F( \tilde \alpha )< 1$ and a second hypothesis ($\E[N^2]<+\infty$) which is satisfied in our work. We use both theorems.
Thanks to the hypothesis of existence of $\psi$ (again by the hypothesis of ellipticity), $\tilde F(\tilde \gamma) \leq 1$ and therefore $\tilde F(\tilde \alpha) \leq 1$. Indeed for all $t>0$ $\tilde F(\tilde \gamma) \leq \E\left[\sum_{|z|=1}\exp(t({V(z)- \tilde \gamma}))\right]$, which by taking the infimum over all $t>0$ in both part of the inequality leads to $\tilde F(\tilde \gamma) \leq \exp(J(\tilde \gamma )) =1$. Moreover  if $\tilde F(\tilde \alpha)>1$, then we should have $\exp{\tilde J(\tilde \gamma)}>1$ which is absurd. \\
Theorem 1 of \cite{McDiarmid}, says that there exists a constant $c_1>0$, such that $P$ almost surely on the set of non-extinction $\max_{|x|=\ell}V(x)-M_\ell \geq  c_1 \log \ell$ with $M_\ell$ the median of $\max_{|x|=\ell}V(x)$, moreover if $\tilde F(\tilde \alpha) = 1$, then  $M_\ell \geq \tilde \alpha \ell-c_1'\log \ell$, with $c_1'>0$. So we only have to check that $\tilde \alpha= \tilde \gamma$. This is an easy computation, indeed we note that 
\begin{align*}
\E\left[\sum_{|z|=1}e^{t(V(z)-\tilde \alpha )}\right]\geq \E\left[\sum_{|z|=1}e^{t(V(z)-\tilde \alpha )}\un_{V(z) \geq \tilde \alpha}\right] \geq \E\left[\sum_{|z|=1}\un_{V(z) \geq \tilde \alpha} \right]=1,
\end{align*}
taking the infimum over all $t>0$, we get $\exp(\tilde J (\tilde \alpha)) \geq 1$ and as $\tilde{J}(a)$ decreases with $a$ and $\tilde {J}(\tilde \gamma)=0$, we get $\tilde \gamma \geq \tilde \alpha$. The other case is pretty similar, let $\epsilon>0$, 
\begin{align*}
\E\left[\sum_{|z|=1}e^{t(V(z)-\tilde \alpha (1+\epsilon) )}\right]=\E\left[\sum_{|z|=1}e^{t(V(z)-\tilde \alpha (1+\epsilon) )}\un_{V(z) \leq \tilde \alpha}\right]+\E\left[\sum_{|z|=1}e^{t(V(z)-\tilde \alpha (1+\epsilon))}\un_{V(z) > \tilde \alpha}\right],
\end{align*}
 as for $\vert z\vert =1$, $V(z)\leq -\log \epsilon_0,$, by definition of $\tilde \alpha$ the last term is equal to $0$, so we get 
 \begin{align*}
\E\left[\sum_{|z|=1}e^{t(V(z)-\tilde \alpha (1+\epsilon) )}\right] \leq  e^{-t\tilde \alpha \epsilon } \E\left[\sum_{|z|=1}\un_{V(z) \leq \tilde \alpha}\right] \leq e^{-t\tilde \alpha \epsilon } e^{\psi(0)}
\end{align*}
taking the infimum over all $t>0$, we get $\exp(\tilde J (\tilde \alpha(1+ \epsilon)))=0$, so $\tilde \gamma \leq \tilde \alpha(1+\epsilon)$.\\
For the case $\tilde F(\alpha) <1$ we use Theorem 2 (b) in \cite{McDiarmid} note that it is the point where we use the hypothesis of second moment for $N$, it gives that there exists a constant $c_2$ such that $P$ almost surely $\max_{|x|=\ell}V(x) \geq \tilde \gamma \ell- c_2 \log n$. \hfill$\blacksquare$


\noindent \\ We finally need a last technical Lemma which tells that, the numbers of back and forth  movement between $z_i$ and $\underline{z}_i$ is small for all $i$.

\begin{Lem} \label{lem2.5} For all the recurrent cases, for all $\epsilon>0$, $\pne$ for $n$ large enough
\begin{align}
\p_\phi^{\mathcal E}\left( \bigcup_{1 \leq j \leq \Ue} \left\{  \sum_{l=1}^{\lo({z}_j,T_{\phi}^{n})-1} \un_{\lo(\underline{z}_j,T_{z_j}^{l+1})-\lo(\underline{z}_j,T_{z_j}^{l})\geq 1} \geq  8/ \epsilon  \right\}  \right) \leq \frac{1}{n^{1+ \epsilon/4}}.
\end{align}
\end{Lem}

\begin{Pre}
Let us denote $ \bigcup_{1 \leq j \leq \Ue} \bar {\mathcal{C}}_j$, the event in the above probability. We have :
\begin{align}
\p_\phi^{\mathcal E}\left( \bigcup_{1 \leq i \leq \Ue} \mathcal \{ \bar {\mathcal{C}}_j,\mathcal A_j \} \right) \leq \sum_{ 1 \leq i \leq \Ue} \p_\phi^{\mathcal E}(Y_{h_n n}(i) \geq 8/ \epsilon )
\end{align}
where $Y_{h_n n}(j):=\sum_{l=1}^{h_n n} \un_{\lo(\underline{z}_j,T_{z_j}^{l+1})-\lo(\underline{z}_j,T_{z_j}^{l})\geq 1}$. By the strong Markov property $Y_{h_n n}(j)$ is a binomial with parameters $h_n n$ and $\tilde{v}_{z_j}:=\p^{\mathcal E}_{z_j}(T_{z_j}>T_{\underline{z}_j})$. As $ \tilde{v}_{z_j} \leq e^{-({V}(\underline{z}_j)-V(\overset{\rightarrow}{z_j}))}$
, so thanks to Lemma \ref{LemVL}, $\pne$ for all $j\leq U_{\epsilon}$ and all $n$ large enough $\tilde{v}_{z_j} \leq  e^{-\log n (1+ \epsilon/4)}$.
Moreover as we have no restriction for $h_n$ but the fact that it goes to infinity with $n$, we can take it for example equal to $\log n$, so we get   that $nh_n   \tilde{v}_{z_j} \leq \log n/n^{\epsilon/4}$. We can now use, for example, the result of \cite{LeCam}, to get that $\pne$ for all $j\leq U_{\epsilon}$ and all $n$ large enough
\begin{align*}
\p_\phi^{\mathcal E}(Y_{h_n n}(j) \geq  8/ \epsilon  ) \leq e^{- \log n/n^{\epsilon/4}}  \left(\frac{\log n}{ n^{\epsilon/4}}\right)^{8/\epsilon}+ 4 \frac{\log n}{n^{1+\epsilon/2}},
\end{align*}
and we conclude with $\p_\phi^{\mathcal E}\left(\cup_{1 \leq j \leq \Ue}  \bar{\mathcal C}_j  \right) \leq \p_\phi^{\mathcal E}\left(\cup_{1 \leq j \leq \Ue} \left \{ \bar{\mathcal C}_j, \mathcal A_j \right \} \right)+ \p_\phi^{\mathcal E}\left(\cup_{1 \leq j \leq \Ue} \bar{\mathcal A}_j \right) $.
\end{Pre}

\noindent \\ Now we move to the proof of the upper bound for $\tilde{R}_n$. Let $\mathcal D_i:= \left\{ \min_{z> \underline{z}_i} \lo(z,T_{\phi}^n) \geq 1\right\}$ such that all $z$ belongs to generation $(1+ \epsilon) \log n$. We have $ \left\{\frac{\tilde{R}_n}{\log n} > c_2\right\} \subset \bigcap_{i=1}^{U_{\epsilon}} \mathcal D_i$. Let us compute an upper bound of the probability $\tpa\left(\cap_{i=1}^{\Ue}\{\mathcal A_i,\mathcal B_i,\mathcal C_i,\mathcal D_i \}\right)$, where $\mathcal A_i$, $\mathcal B_i$, and $\mathcal C_i$ have been defined in the previous Lemmata. 
 We have 
\begin{align*}
 &\tpa\left(\bigcap_{i=1}^{\Ue}\{\mathcal A_i,\mathcal B_i,\mathcal C_i,\mathcal D_i \}\right)=
 \prod_{j=1}^{\Ue} \sum_{k_j=1}^{h_n n-1}\sum_{l_j=0}^{  8 / \epsilon-1}\sum_{m_j=l_j}^{K_n h_n n-1}\\ 
 &\tpa\left(\bigcap_{i=1}^{\Ue} \left\{\mathcal D_i, \lo(z_i,T_{\phi}^n)=k_i, \lo(\underline{z}_i,T_{z_i}^{k_i})=m_i,\sum_{l=1}^{k_i-2} \un_{\lo(\underline{z}_i,T_{z_i}^{l+1})\neq\lo(\underline{z}_i,T_{z_i}^{l})} =l_
i \right\}  \right) .
\end{align*} 
In the following expression we add a sum over all the possible sequences  $(q_1^i,\cdots,q_{l_i}^i)$ of the different time of excursions from $\underline{z}_i$ to $z_i$: for this we denote $\mathcal G_i^{m_i}(q_1^i,\cdots,q_{l_i}^i)$ the event that says that during the $m_i$ returns to $\underline{z}_i$, the walk will touch the point ${z}_i$ only between the $(q_r^i-1)$nth and $q_r^i${nth} return time to $\underline{z}_i$ for all $r \leq l_i$.   
 \begin{align*}
& \tpa\left(\bigcap_{i=1}^{\Ue}\{\mathcal A_i,\mathcal B_i,\mathcal C_i,\mathcal D_i \}\right)\\
& = \prod_{j=1}^{\Ue} \sum_{k_j=1}^{h_n n-1}\sum_{l_j=0}^{ 8 / \epsilon-1}\sum_{m_j=l_j}^{K_n h_n n-1} \sum_{q_1^j,\cdots,q_{l_j}^j}\tpa\left(\bigcap_{i=1}^{\Ue} \left\{\mathcal D_i, \lo(z_i,T_{\phi}^n)=k_i, \lo(\underline{z}_i,T_{z_i}^{k_i})=m_i,\mathcal G_i^{m_i}(q_1^i,\cdots,q_{l_i}^i) \right\}  \right). 
\end{align*} 
 Now on $\left\lbrace \mathcal G_i^{m_i}(q_1^i,\cdots,q_{l_i}^i), \lo(z_i,T_{\phi}^n)=k_i,\lo(\underline{z}_i,T_{z_i}^{k_j})=m_i\right\rbrace$ the event $\mathcal D_i$ can be written 
  \begin{align*}
  \mathcal D_i&=\left\{ \min_{z> \underline{z}_i} \left(\sum_{s_i=0}^{l_i-1} \lo(z,T_{\underline{z}_i}^{q_{s_i+1}^i-1})-\lo(z,T_{\underline{z}_i}^{q_{s_i}^i}) \right) \geq 1\right\} \\
  & = \bigcup_{s_i=0}^{l_i-1} \left\{\min_{z> \underline{z}_i}  \left(  \lo(z,T_{\underline{z}_i}^{q_{s_i+1}^i-1})-\lo(z,T_{\underline{z}_i}^{q_{s_i}^i}) \right) \geq 1\right\}  =: \bigcup_{s_i=0}^{l_i-1} \mathcal H_i(q_{s_i}^i) .
  \end{align*} 
 We finally get  
  \begin{align*}
&\tpa\left(\bigcap_{i=1}^{\Ue} \left\{\mathcal D_i, \lo(z_i,T_{\phi}^n)=k_i, \lo(\underline{z}_i,T_{z_i}^{k_i})=m_i,\mathcal G_i^{m_i}(q_1^i,\cdots,q_{l_i}^i) \right\}  \right)\\
& \leq \prod_{j=1}^{\Ue}\sum_{s_j=0}^{l_j-1} \tpa\left(\bigcap_{i=1}^{\Ue} \left\{\mathcal H(q_{s_i}^i), \lo(z_i,T_{\phi}^n)=k_i, \lo(\underline{z}_i,T_{z_i}^{k_i})=m_i,\mathcal G_i^{m_i}(q_1^i,\cdots,q_{l_i}^i) \right\}  \right) \\
 & \leq \prod_{j=1}^{\Ue}\sum_{s_j=0}^{l_j-1} \tpa\left(\bigcap_{i=1}^{\Ue} \left\{\mathcal H(q_{s_i}^i), \lo(z_i,T_{\phi}^n)=k_i, \lo(\underline{z}_i,T_{z_i}^{k_i})=m_i,\tilde{\mathcal G_i} (q_{s_i}^i) \right\}  \right) 
\end{align*} 
where  $\tilde{\mathcal G_i} (q_{s_i}^i):= \{\forall r, T_{\underline{z}_i}^{q_{s_i+1}^i-1} \leq r \leq T_{\underline{z}_i}^{q_{s_i}^i},\ X_r> z_i \}$. The next step is to make disappear $\lo(z_i,T_{\phi}^n)=k_i$, and $\lo(\underline{z}_i,T_{z_i}^{k_i})=m_i$ carefully, we simply notice that  
 \begin{align*}
& \prod_{j=1}^{\Ue} \sum_{k_j=1}^{h_n n}\tpa\left(\bigcap_{i=1}^{\Ue} \left\{\mathcal H(q_{s_i}^i), \lo(z_i,T_{\phi}^n)=k_i, \lo(\underline{z}_i,T_{z_i}^{k_i})=m_i,\tilde{\mathcal G_i} (q_{s_i}) \right\}  \right) \\
& \leq  \prod_{j=1}^{\Ue} \sum_{k_j=1}^{h_n n} \tpa\left(\bigcap_{i=1}^{\Ue} \left\{\mathcal H(q_{s_i}^i), \lo(z_i,T_{\phi}^n)=k_i,\tilde{\mathcal G_i} (q_{s_i}) \right\}  \right) 
=\tpa\left(\bigcap_{i=1}^{\Ue} \left\{\mathcal H(q_{s_i}^i),\tilde{\mathcal G_i} (q_{s_i}) \right\}  \right). 
\end{align*} 
We are now ready to apply the strong Markov property, indeed the $ (T_{\underline{z}_i}^{q_{s}^i},i \leq \Ue)$ can now be ordered, and as they are stopping times recursively we finally get:
\begin{align*}
 \tpa\left(\bigcap_{i=1}^{\Ue} \left\{\mathcal H(q_{s_i}^i),\tilde{\mathcal G_i} (q_{s_i}) \right\}  \right)&=\prod_{i=1}^{\Ue} \tpa_{\underline{z}_i}\left(\min_{z> \underline{z}_i}    \lo\left(z,T_{\underline{z}_i}^{q_{s_i+1}^i-1-q_{s_i}^i} \right) \geq 1  \right)  \\
 & \leq \prod_{i=1}^{\Ue} \tpa_{\underline{z}_i}\left(\min_{z> \underline{z}_i}    \lo\left(z,T_{\underline{z}_i}^{K_n h_n n} \right) \geq 1  \right)  .
\end{align*} 
We are left to get an upper bound for the probabilities in the above product, and also to count the number of term we have in the previous product of sums. First about the sums we notice that $\sum_{l_1=0}^{8/\epsilon}\sum_{m_1=l_1}^{K_n h_n n}\sum_{q_1^1,\cdots,q_{l_1}^1} \sum_{s_1=0}^{l_1}1=\sum_{l_1=0}^{8/\epsilon}\sum_{m_1=l_1}^{K_n h_n n}\binom{m_1}{l_1} (l_1+1) \leq (8/\epsilon+1)K_n h_n n \sum_{l_1=0}^{8/\epsilon}\binom{K_n h_n n }{l_1} \leq (8/\epsilon+1)^2 K_n h_n n (K_n h_n n)^{8/\epsilon}$, so finally
$$
 \prod_{j=1}^{\Ue} \sum_{l_j=0}^{  8 / \epsilon}\sum_{m_j=l_j}^{K_n h_n n} \sum_{q_1^j,\cdots,q_{l_j}^j}\sum_{s_j=0}^{l_j}1 \leq \left((8/\epsilon+1)^2 K_n h_n n (K_n h_n n)^{8/\epsilon} \right)^{\Ue} \sim (K_n h_n n)^{\Ue \left(8/\epsilon+1\right)}.$$
Using successively the strong Markov property,  (\ref{eq2}) and  the hypothesis of ellipticity for all $z>z_i$:
\begin{eqnarray*}
\tpa_{\underline{z}_i}\left(\min_{z> \underline{z}_i}    \lo\left(z,T_{\underline{z}_i}^{K_n h_n n} \right) \geq 1  \right)  &\leq&  \tpa_{\underline{z}_i}\left(  \lo\left(z,T_{\underline{z}_i}^{K_n h_n n} \right) \geq 1  \right)
=1- \tpa_{\underline{z}_i}\left(T_{\underline{z}_i}<T_z\right)^{K_n h_n n}\\
&Ê\leq &1-  \left(1-p(\underline{z}_i,\underline{\qz}_i)\frac{1}{\sum_{u\in\rrbracket \underline{z}_i,z\rrbracket }e^{V(u)-V(\underline{\qz}_i)}}\right)^{K_n h_n n} \\
& \leq & 1- \exp\left(-c K_n h_n n  e^{-\max_{u\in\rrbracket \underline{z}_i,z\rrbracket}V(u)-V(\underline{\qz}_i)} \right),
\end{eqnarray*}
with $c>0$. The stationarity of $V$ gives the following equality in law with respect to $P$:
$\max_{u\in\rrbracket \underline{z}_i,z\rrbracket}V(u)-V(\underline{\qz}_i)=\max_{|z|=\frac{\epsilon}{2\tilde \gamma} \log n }V(z)$, moreover thanks to lemma \ref{LemVL}, $\pne$ for all $n$ large enough: 
$$\max_{|z|=\frac{\epsilon}{2\tilde \gamma} \log n}V(z) \geq (1-\epsilon)\frac{\epsilon}{2} \log n.$$
We finally get that $\pne$ for all $n$ large enough: 
\begin{eqnarray*}
\tpa_{\underline{z}_i}\left(\min_{z> \underline{z}_i}    \lo\left(z,T_{\underline{z}_i}^{K_n h_n n} \right) \geq 1  \right)  &\leq& \frac{c' K_n h_n n}{n^{\epsilon(1-\epsilon)/2}},
\end{eqnarray*}
implying
\begin{align*}
 \tpa\left(\bigcap_{i=1}^{\Ue} \left\{\mathcal H(q_{s_i}^i),\tilde{\mathcal G_i} (q_{s_i}) \right\}  \right)& \leq  \left(\frac{c' K_n h_n n}{n^{\epsilon(1-\epsilon)/2}}\right)^{\Ue}. 
 \end{align*}
Collecting all what we did above and replacing $K_n h_n n$ by its value, yields that $\pne$ for $n$ large enough
\begin{align*}
 \tpa\left(\bigcap_{i=1}^{\Ue}\{\mathcal A_i,\mathcal B_i,\mathcal C_i,\mathcal D_i \}\right) \leq (K_n h_n n)^{\Ue(8/\epsilon+2)} \left(\frac{1}{n^{\epsilon(1-\epsilon)/2}}\right)^{\Ue}.
 \end{align*}
From Kesten-Stigum theorem \cite{KesSti} (here the hypothesis that $E(N \log^+ N)<\infty$ is trivially satisfied), we know that $\pne$  $\lim_{\epsilon \rightarrow 0}\Ue/e^{\psi(0)u_{\epsilon}} =W$ where $W$ a strictly positive, finite random variable. In particular choosing $u_{\epsilon}=\frac{1}{e^{\psi(0)}}\log \frac{1}{\epsilon^2}$ $\pne$ for all $\epsilon>0$ small enough  $ 4/ (1-\epsilon) \epsilon \leq \Ue \leq 1/\epsilon^3 $, finally  remember that $K_n$ is given just after \ref{alakokan}  so we get $\pne$ for $n$ large enough 
 \begin{align*}
 \tpa\left(\bigcap_{i=1}^{\Ue}\{\mathcal A_i,\mathcal B_i,\mathcal C_i,\mathcal D_i \}\right) \leq \frac{(\log n)^{c^{\prime\prime} \epsilon^4}}{n^{2}}, 
 \end{align*}
with $c''>0$. Finally collecting the result of the different Lemmata we get that $\pne,$  $ \tpa\left(\frac{\tilde{R}_{n}}{\log n} > c_2\right)$ is summable, applying Borel-Cantelli Lemma we get \ref{limsup}.

\section{Connexion between  \boldmath{$\tilde{R}_n$} and  \boldmath{$R_n$}}

\subsection{Case  \boldmath{$\psi(1)=0,\ \psi'(1)\geq 0$} or  \boldmath{$\inf_{t \in [0,1]}\psi(t)<0$}}

We have the following

\begin{Lem} \label{lem1.1} Assume \ref{hypP} or \ref{hypPP} or \ref{hyp0}, then for all $\epsilon>0$ $\pn$ for all $n$ large enough
\begin{eqnarray}
 \tilde{R}_{n^{1- \epsilon}} \leq R_n \leq \tilde{R}_n.
\end{eqnarray}
\end{Lem}
Note that only the first inequality needs to be proven, moreover the case (\ref{hypP}) and (\ref{hypPP}) follows directly by the fact that the random walks are positive recurrent. In what follows we will always assume that \ref{hyp0} is realised and for $m\in\mathbb N$, we denote $ \mathcal{T}_m:=\inf\lbrace k\geq 0, \vert X_k\vert =m \rbrace$ the hitting time of the generation $m$. The key-point is the following
\begin{Lem}\label{lemfin} There exists a constant $\alpha>0$, such that  $\pn$ for all $m$ large enough
\begin{eqnarray}
\mathcal{A}_1(m)&:=& \{ \lo(\phi,\mathcal{T}_m) \geq \exp{((m\alpha)^{1/3}(1-\epsilon/4))} \}, \label{1.2} and \\
\mathcal{A}_2(m)&:=&\{\mathcal{T}_m \leq \exp{((m\alpha)^{1/3}(1+\epsilon/2))}\} \label{1.3}
\end{eqnarray}
are realized.
\end{Lem}

\noindent From the above Lemma the proof of the first Lemma is straightforward, indeed for $n$ large enough on $\mathcal{A}_2$ 
\begin{align*}
  \lo(\phi,n) \geq \lo\left(\phi, \mathcal{T}_{\frac{(\log n)^3}{ \alpha(1+\nicefrac{\epsilon}{2})^3}}\right),
\end{align*}
therefore, for $n$ large enough on $\mathcal{A}_1$ and $\mathcal{A}_2$ 
\begin{align*}
  \lo(\phi,n) \geq \exp( \log n(1-\epsilon/4)/(1+\epsilon/2))\geq n^{1-\epsilon}.
\end{align*}

So we are left to prove Lemma \ref{1.2}, notice that it can be deduced from what is done in \cite{HuShi10b}, for completness we give some details here except the proof of the following delicate to prove Lemma 

\begin{Lem} \label{lemmil} (\cite{HuShi10b}) For all $\epsilon>0$ $\pne$ for all $m$ large enough
 \begin{eqnarray}
 \rho_m \leq \exp(-m^{1/3} \alpha^{1/3}(1-\epsilon/8 )),
\end{eqnarray} 
where $ \rho_m:=\p_{\phi}(\mathcal{T}_m<T_{\phi})$. 
\end{Lem}

\noindent\textbf{Proof of lemma \ref{lem1.1}} For $\mathcal{A}_1(m)$, the strong Markov property gives $\tpa\left(\lo(\phi,\mathcal{T}_m)\geq k\right)= (1- \rho_m)^k $, then Lemma \ref{lemmil} yields that $\pne$ for $m$ large enough
\begin{eqnarray}
\tpa\left(\lo(\phi,\mathcal{T}_m)\leq \exp(m^{1/3} \alpha^{1/3}(1-\epsilon/4 )) \right) 
& \leq & \exp(-m^{1/3} \alpha^{1/3} \epsilon/8) ,
\end{eqnarray}
applying Borel-Cantelli Lemma leads to \ref{1.2}.  \\
For $\mathcal{A}_2(m)$, from U.A. Rozikov \cite{Rozikov}, $\E^{\mathcal E}[\mathcal{T}_m]=\frac{\gamma_m(\phi)}{\rho_m}$, where $\gamma_m(\phi)$ is defined in the appendix. Lemma \ref{lema3} and \ref{lemmil} imply the existence of a constant $c'>0$ such that $\pne$ for $m$ large enough
\begin{eqnarray}
\E^{\mathcal E}[\mathcal{T}_m] \leq c' m\exp( m^{1/3} \alpha^{1/3}(1+\epsilon/8 )),  
\end{eqnarray}
the Markov  inequality together with the above inequality yields that $\pne$ $\tpa(\mathcal{T}_m > m^{1/3} \alpha^{1/3}(1+\epsilon/4 ))$ is summable and we conclude with Borel-Cantelli Lemma.\hfill
$\blacksquare$

\noindent \\
Finally notice that by Lemma \ref{lem1.1}, (\ref{liminf}) and $(\ref{limsup})$, $\pn$ for $n$ large enough
\begin{align}
 \frac{1}{\tilde \gamma} (1-\epsilon)^{2}\leq \frac{R_n}{\log n} \leq \frac{1}{\tilde \gamma}(1+\epsilon)
\end{align}
we get the Theorem for the first three cases by letting $\epsilon$ go to zero.


\subsection{Case \boldmath{$\psi(1)=0,\ \psi'(1) < 0$} }

Let us prove
\begin{Lem} \label{lem10} Under $\psi(1)=0,\ \psi'(1) < 0$, we have for all $\epsilon>0$, $\pn$ for all $n$ large enough
\begin{eqnarray*}
 \tilde{R}_{n^{\nu'-\epsilon}} \leq R_n \leq \tilde{R}_{n^{\nu'+\epsilon}},
\end{eqnarray*}
where $\nu':=1/\min(\kappa,2)$.
\end{Lem}
To prove this Lemma we use the following results of \cite{HuShi10} that can be extended to a supercritical Galton Watson tree by using the same technics:  
\begin{Pro} \label{Prop1} (\cite{HuShi10}) Under $\psi(1)=0,\ \psi'(1) < 0$, we have for all $\epsilon>0$, $\pn$ for all $m$ large enough
\begin{align}
& m^{-\epsilon} E[\beta_m(\phi^{1})] \leq   \beta_m(\phi^{1}) \leq m^{\epsilon}E[\beta_m(\phi^{1})], 
\end{align}
where $\beta_m(\phi^{1}):=\tpa_{\phi^{1}}[\mathcal{T}_m \leq T_{\phi}]$. Moreover if $\kappa \in (2,+ \infty]$, $ E[\beta_m(\phi^{1})] = O(1/m)$ and if $\kappa \in (1,2]$ $m^{- \frac{1}{\kappa-1}-\epsilon}\leq E[\beta_m(\phi^{1})] \leq m^{- \frac{1}{\kappa-1}
+\epsilon}$.
Also $\pn$ for all $n$ large enough
\begin{align}
& \mathcal{T}_{n^{\nu(1-\epsilon)}}\leq n \leq \mathcal{T}_{n^{\nu(1+\epsilon)}},
\end{align}
with $\nu:=1-1/ \min\{\kappa,2\}$.
\end{Pro}

\textbf{Proof of lemma \ref{lem10}} First notice that thanks to the second part of the above proposition,  $\pn$ for all $n$ large enough
\begin{align}
 \lo(\phi,\mathcal{T}_{n^{\nu(1-\epsilon)}}) \leq \lo(\phi,n) \leq \lo(\phi,\mathcal{T}_{n^{\nu(1+\epsilon)}}). \label{3.23}
\end{align}
\textit{The upper bound}
we study the asymptotic of $ \lo(\phi,\mathcal{T}_{m})$ for large $m$, using Markov inequality we have 
\begin{align}
\tpa\left(\lo(\phi,\mathcal{T}_{m})\geq m^{2 \epsilon}/E[\beta_m(\phi^{1})] \right) \leq \frac{(1-\rho_m)E[\beta_m(\phi^{1})]}{\rho_mm^{2 \epsilon}}.
\end{align}
By definition $\rho_m=\sum_{i=1}^{N^{(\phi)}}p(\phi,\phi^{(i)})\beta_m(\phi^{i})$, then by using the fact that the $\beta_m(\phi^{(i)})$ are i.d. with mean $E[\beta_m(\phi^{1})]$, the hypothesis of ellipticity  and the first part of the above Proposition we get that there exist positive constants $c_1>0$ and $c_2>0$ such that $\pne$ $ c_1 E[\beta_m(\phi^{1})] m^{-\epsilon} \leq  \rho_m \leq c_2 E[\beta_m(\phi^{1})] m^{\epsilon}$,  so $\pne$ for $n$ large enough
\begin{align}
p_m:=\tpa\left(\lo(\phi,\mathcal{T}_{m})\geq \frac{m^{2 \epsilon}}{ E[\beta_m(\phi^{1})]} \right) \leq  \frac{1}{m^{\epsilon}}. \label{rho1}
\end{align}
We deduce from that the convergence of the sum $\sum_{\ell}p_{\ell^{2/\epsilon}}$, therefore according to Borel Cantelli Lemma $\pn$ for all $l$ large enough $\lo(\phi,\mathcal{T}_{\ell^{\frac{2}{\epsilon} }})\leq {\ell^{\frac{2}{\epsilon}2 \epsilon}}/{ E\left[\beta_{\ell^{\frac{2}{\epsilon}}}(\phi^{1})\right]}$. Taking $ (\ell-1)^{2/\epsilon} \leq m \leq \ell^{2/\epsilon} $ in such a way that for $\ell$ large enough $ \ell^{2/\epsilon} \leq m^{1+\epsilon} $, we get by using the fact that $\lo(\phi,T_\ell)$ is increasing in $\ell$ and $\beta_\ell$ decreasing in $\ell$,  that $\pn$ for all $m$ large enough $\lo(\phi,\mathcal{T}_{m}) \leq {m^{3 \epsilon}}/{E[\beta_{m^{1+\epsilon}}(\phi^{1})]}$. Finally $\pn$  for all $n$ large enough $ \lo(\phi,\mathcal{T}_{n^{\nu(1+\epsilon)}}) \leq {n^{4\nu \epsilon}}/{E[\beta_{n^{(1+3\epsilon)\nu}}(\phi^{1})]}$. Now, distinguishing the two cases we get for  $\kappa \in (1,2]$, $\pn$ for $n$ large enough, $ \lo(\phi,\mathcal{T}_{n^{\nu(1+\epsilon)}}) \leq {n^{\frac{1}{\kappa}+c_0\epsilon}}$, and  for $\kappa \in (2,+ \infty]$, $ \lo(\phi,\mathcal{T}_{n^{\nu(1+\epsilon)}}) \leq {n^{\frac{1}{2}+c_0'\epsilon}}$ where $c_0$ and $c_0'$ are two positive constant. Collecting this result and the right-hand side of \ref{3.23} gives the upper bound.\\
\textit{The lower bound}, let $(\lambda_m,m)$ a positive sequence decreasing to zero when $m$ goes to infinity. First notice that 
\begin{align*}
\E^{\mathcal E}\left[e^{-\lambda_m \lo \left(\phi,\mathcal{T}_m \right)}\right]= \frac{\rho_m}{1-e^{-\lambda_m}(1-\rho_m)},
\end{align*}
therefore for $m$ large enough and by taking $\lambda_m=m^{\epsilon} \rho_m$ we get $\E^{\mathcal E}\left[e^{-\lambda_m \lo \left(\phi,\mathcal{T}_m \right)}\right]\leq  { 2 \rho_m}/{(\lambda_m+\rho_m)} \leq 2m^{-\epsilon} $. We obtain that $\E^{\mathcal E}\left[ \sum_\ell e^{-\lambda_{m_\ell} \lo \left(\phi,\mathcal{T}_{m_\ell} \right)}\right]$ is finite, for the subsequence $m_\ell=\lfloor \ell^{2/\epsilon}\rfloor$,  therefore $\pn$ for all $\ell$ large enough $\lambda_{m_\ell} \lo\left(\phi,\mathcal{T}_{m_\ell} \right)\geq 1$, then it is clear that for all $m\in [m_\ell,m_{\ell+1}]$, $ \lambda_{m_\ell}  \lo \left(\phi,\mathcal{T}_{m} \right)\geq  1$. Moreover using the estimates of $\rho_.$ just above (\ref{rho1}) and of $E[\beta_m(\phi^{1})] $, $\pne$  for all $\ell$ large enough and for all $m\in [m_\ell,m_{\ell+1}]$, $1/ \lambda_{m_\ell} \geq  1/(m^{c_3 \epsilon}\lambda_m)$, with $c_3>0$ a well chosen constant. Therefore for some positive constant $c_4$, $\pn$ for $n$ large enough $ \lo\left(\phi,\mathcal{T}_{m} \right) \geq  \frac{1}{m^{c_4\epsilon}} \frac{1}{E[\beta_{m}(\phi^{1})]}  $. Then we separate the two cases and  use the left hand side of (\ref{3.23}) to get the lower bound. \hfill$\blacksquare$ \\
Lemma \ref{lem10} together with Proposition \ref{Pro1} yields the theorem for this last case.\\
Finally note that Proposition \ref{profin} is a simple consequence of Lemma \ref{lemfin} and proof of lemma  \ref{lem10}.

\section{Appendix}
In this appendix, for completness, we describe and sketch the proof of some classical results. Given a vertex $x\in\mathbb T,$  we denote $x_0:=\phi,\dots,x_n:=x$ the vertices on $\llbracket \phi,x \rrbracket$ with $\vert x_i\vert=i$ for all $0\leq i\leq n$. 

\subsection{Biggins-Kyprianou identities}
For any $n\geq1$ and any mesurable function $F:\mathbb R^n\times\mathbb R^n\rightarrow[0,+\infty)$, Biggins-Kyprianou identity is given by
\begin{equation}\label{B1}
E\left[\sum_{\vert x\vert =n}e^{-V(x)-\psi(1)n}F(V(x_i),1\leq i\leq n)\right]=E[F(S_i,1\leq i\leq n)]
\end{equation}
where $(S_i-S_{i-1})_{i\geq 1}$, are i.i.d. random vectors, and the distribution of $S_1$ is determined by :
\begin{equation}\label{B2}
E[f(S_1)]=E\left[\sum_{\vert x\vert =1}e^{-V(x)-\psi(1)}f(V(x))\right], 
\end{equation}
for any measurable function $f:\mathbb R\rightarrow[0,+\infty)$. A proof can be found in \cite{BigginsKyprianou}, see also \cite{Shi4}.

\subsection{Classical results about birth and  death chains}
\begin{Lem}\label{lemfinal}
For $x^\prime\in\llbracket \phi,x\rrbracket$:
\begin{eqnarray}
\p_{x^\prime_x}^{\mathcal E}(T_x<T_{x^\prime})&=&\frac{e^{V(x^\prime_x)}}{\sum_{z\in\rrbracket x^\prime,x\rrbracket }e^{V(z)}} \label{eq1},\\
\p_{\overset{\leftarrow}{x}}^{\mathcal E}(T_{x^\prime}<T_x)&=&\frac{e^{V(x)}}{\sum_{z\in\rrbracket x^\prime,x\rrbracket }e^{V(z)}}.\label{eq2}
\end{eqnarray}
where $x^\prime_x$ is the only children of $x^\prime$ in $\llbracket x^\prime,x\rrbracket$. 
\end{Lem}

\textbf{Proof:}
Let $(\sigma_n)_{n\geq0}$ the family of stopping times defined by $\sigma_n=\inf\lbrace k>\sigma_{n-1},X_k\in\llbracket \phi,x\rrbracket,X_k\neq X_{{\sigma_{n-1}}}\rbrace $ and define $Z_n=X_{\sigma_n}$ for $n\geq0$.
$(Z_n)_{n\geq 0}$ is a birth and death Markov chain on $\llbracket \phi, x\rrbracket$ with transition probabilities given by:
\begin{eqnarray*}
p_{x_i}&:=&\p^{\mathcal E}(Z_{n+1}=x_{i+1}\vert Z_n=x_i)=\frac{A(x_{i+1})}{1+A(x_{i+1})},\\
q_{x_{i}}&:=&\p^{\mathcal E}(Z_{n+1}=x_{i-1}\vert Z_n=x_i)=\frac{1}{1+A(x_{i+1})},
\end{eqnarray*}
$\forall 1\leq i\leq n-1$ and $p_{\phi}=q_{x}=1$,
indeed
\begin{eqnarray*}
p_{x_i}
&=&\p_{x_i}^{\mathcal E}(X_{\sigma_1}=x_{i+1})
=\sum_{\ell\geq0}\p_{x_i}^{\mathcal E}(X_{T_{x_i}^\ell+1}=x_{i+1}, \forall m<\ell, X_{T_{x_i}^m+1}\notin\llbracket\phi,x\rrbracket)\\
&=&\sum_{\ell\geq0}p({x_i,x_{i+1}})\p_{x_i}^{\mathcal E}(X_1\notin\llbracket\phi,x\rrbracket)^\ell=\frac{p({x_i,x_{i+1}})}{1-\p_{x_i}^{\mathcal E}(X_1\notin\llbracket\phi,x\rrbracket)}=\frac{p({x_i,x_{i+1}})}{1-\sum_{k\neq j}p({x_i,x_i^{(k)}})}\\
&=&\frac{A(x_{i+1})}{1+A(x_{i+1})}.
\end{eqnarray*}
Let us introduce:
\begin{equation*}
\xi_{0}:=1,\,\xi_{\ell}:=\prod_{k=1}^{\ell}\frac{q_{k}}{p_{k}},
\,\ell\geq 1,
\end{equation*}
and consider $f:\mathbb N\rightarrow\mathbb R$ given by $f(\phi)=0$ and for $1\leq k\leq n,f(x_k)=\sum_{\ell=0}^{k-1}\xi_\ell$. Easily we can see that $(f(Z_k))_{k\geq0}$ 
is a martingale.  With $\tau_i=\inf\lbrace m\geq, 0, Z_m=x_i\rbrace$ and for $1\leq i<j<k$, according to the optional stopping time Theorem, for $1\leq i<j<k$ :
\begin{eqnarray*}
f(x_j)
&=&\E^{\mathcal E}_{x_j}\left[f(X_{\tau_{i}\wedge \tau_{k}})\right]
=f(x_i)\p^{\mathcal E}_{x_j}(\tau_{i}<\tau_{k})+f(x_k)[1-\p^{\mathcal E}_{x_j}(\tau_{i}<\tau_{k})]\\
&\Leftrightarrow&\p^{\mathcal E}_{x_j}(\tau_{i}<\tau_{k})
=\frac{\sum_{\ell=i}^{j-1}\xi_{\ell}}{\sum_{\ell=i}^{k-1}\xi_{\ell}}=\frac{\sum_{z\in\rrbracket x_i,x_j\rrbracket}e^{V(z)}}{\sum_{z\in\rrbracket x_i,x_k\rrbracket}e^{V(z)}}
\end{eqnarray*}
recalling that $V(x)=-\sum_{z\in\rrbracket \phi,x\rrbracket }\log A(z),x\in\mathbb T\backslash \lbrace\emptyset\rbrace$.
Since $\lbrace\tau_{x}<\tau_{x^\prime}\rbrace=\lbrace T_{x}<T_{{x^\prime}}\rbrace$ conditionnaly on $\lbrace X_0={x^\prime}_x\rbrace$, thus formula \ref{eq1} is proved.
\hfill$\blacksquare$\\

\subsection{About  \boldmath{$(\gamma_n,n)$}}
Let us define:
\begin{equation}\label{gamma1}
 \gamma_n(x):=\left\lbrace\begin{array}{cl}
0&\mbox{ if $\vert x\vert=n$, }\\
 \frac{ \nicefrac{1}{p(x,\overset{\leftarrow }{x})}+\sum_{i=1}^{N_x}A(x^i)\gamma_n(x^i)}{1+\sum_{i=1}^{N_x}A(x^i)\beta_n(x^i)},&\mbox{ if $1\leq \vert x\vert<n$,  }Ê\\
\sum_{i=1}^Np(\phi,\phi_i)\gamma_n(\phi_i), &\mbox{ if $x=\phi$ }Ê.
\end{array}\right.
\end{equation}
where $\beta_n:=\p_x^{\mathcal E}(\mathcal T_n<T_{\overset{{\leftarrow}}{x}})$.
\begin{Lem} \label{lema3} 
Assuming $\psi(1)=0$:
\begin{eqnarray}\label{gamman}
& & \sup_{n \geq 1} \frac{\gamma_n(\phi)}{n}< + \infty, \p.a.s.
\end{eqnarray}
\end{Lem}

This result is already proved in the case of a $b$-ary tree (see for instance  \cite{HuShi10a}). Here, we treat the case of a Galton-Watson tree.\\
\textbf{Proof:}

First, observe that for all $2\leq k\leq n$ :
\begin{equation}\label{rec}
\gamma_n(\phi)\leq K\sum_{j=1}^{k-1}\sum_{\vert x\vert =j}\prod_{y\in\rrbracket \phi  ;x \rrbracket}A(y)+\sum_{\vert x\vert =k}\left(\prod_{y\in\rrbracket \phi  ;x \rrbracket}A(y)\right)\gamma_n(x)
\end{equation} 
where $K$ is a constant satisfying $\forall x\in \mathbb{T},p(x,\overset{\leftarrow}{x})^{-1}\leq K$. The existence of $ K$ is provided by assumptions \ref{hyp1}.\\
As $p(\phi,\phi^{i})\leq A(\phi^{i}),\,\forall 1\leq i\leq N $, we deduce from (\ref{gamma1}):
\begin{equation}\label{ineg1}
\gamma_n(\phi)\leq\sum_{i=1}^NA(\phi^{i})\gamma_n(\phi^{i}), 
\end{equation}
and note that formula (\ref{gamma1}) implies :
\begin{equation}\label{ineg2}
\gamma_n(x)\leq K+\sum_{i=1}^{N_x}A(x^{i})\gamma_n(x^{i}),   \forall 1\leq \vert x\vert\leq  n.
\end{equation}
Then from (\ref{ineg1}) and (\ref{ineg2}), we deduce formula (\ref{rec}) for $k=2$:
\begin{eqnarray*}
\gamma_n(\phi)
&\leq&\sum_{i=1}^NA(\phi^{i})(K+\sum_{j=1}^{N_{\phi^{i}}}A({\phi^{i,j}})\gamma_n({\phi^{i,j}})= K\sum_{i=1}^NA(\phi^{i})+\sum_{i=1}^N\sum_{j=1}^{N_{\phi^{i}}}A(\phi^{i})A({\phi^{i,j}})\gamma_n({\phi^{i,j}})\\
&=&K\sum_{\vert x\vert =1}\prod_{y\in\rrbracket \phi  ;x \rrbracket}A(y)+\sum_{\vert x\vert =2}\left(\prod_{y\in\rrbracket \phi  ;x \rrbracket}A(y)\right)\gamma_n(x)
\end{eqnarray*}
Assume that (\ref{rec}) is true for  one $k\geq2$ 
, we prove that it still true for $k+1$. Using again (\ref{ineg2}):
\begin{eqnarray*}
\gamma_n(\phi)&\leq &K\sum_{j=1}^{k-1}\sum_{\vert x\vert =j}\prod_{y\in\rrbracket \phi  ;x \rrbracket}A(y)+\sum_{\vert x\vert =k}\left(\prod_{y\in\rrbracket \phi  ;x \rrbracket}A(y)\right)\left(K+\sum_{i=1}^{N_x}A(x^{i})\gamma_n(x^{i})\right)\\
&\leq&K\sum_{j=1}^{k}\sum_{\vert x\vert =j}\prod_{y\in\rrbracket \phi  ;x \rrbracket}A(y)+\sum_{\vert x\vert =k+1}\left(\prod_{y\in\rrbracket \phi  ;x \rrbracket}A(y)\right)\gamma_n(x)
\end{eqnarray*}

Applying formula (\ref{rec}) to $k=n$ and recalling that $\gamma_n(x)=0$ for $\vert x\vert=n$ :
\begin{equation}
\gamma_n(\phi)\leq K\sum_{j=1}^{n-1}\sum_{\vert x\vert =j}\prod_{y\in\rrbracket \phi  ;x \rrbracket}A(y)=K\sum_{j=1}^{n-1}M_j,\label{ineg3}
\end{equation}
where $M_j:=\sum_{\vert x\vert=j}\prod_{\rrbracket \phi; x\rrbracket}A(y)$.  $(M_j)_{ j\geq1} $ is a positive $\mathcal F_j$-martingale with $M_0=1$ and $\mathcal F_j:=\sigma\lbrace (A(x^{1}),\cdots, A(x^{N_x}),N_x): \vert x\vert \leq j,x\in \mathbb{T}\rbrace$:
\begin{itemize}
\item obviously we have positivity and for all $j\geq0,\, M_j\in \mathcal F_j$;
\item for all  $x\in T$, as $(A(x^{1}),\cdots, A(x^{N_x}),N_x)$ is equal in law to the vector $(A_1,\cdots,A_N,N)$:
\begin{eqnarray*}
\E[M_{j+1}\vert \mathcal F_j]
=M_j\E[\sum_{i=1}^{N}A_i],
\end{eqnarray*}
and we conclude with $M_0=\E[\sum_{i=1}^{N}A_i]=1$, since $\psi(1)=0$.
\end{itemize}
Consequently, there exists an almost sure limit for $(M_j)_{j\geq0}$ which implies that  $\sup_{j}M_j<\infty$ almost surely.\\
 Thus, (\ref{ineg3}) implying  $\frac{\gamma_n(\phi)}{n}\leq K \sup_{j}M_j$, the proof is complete.
$\square$\\

\vspace{0.5cm}
\noindent \textbf{Acknowledgments}
We would like to thank J-B Gou\'er\'e for sharing several discussions on  branching random walks, and Thomas Haberkorn for nice numerical simulations.

 \bibliography{thbiblio}
\end{document}